\documentclass[a4paper]{amsart} 

\usepackage{amssymb} \usepackage{amscd} \usepackage{times}

\numberwithin{equation}{subsection}

\def\zbb{\mathbb{Z}}  
  
  \def\phi{\varphi}
 \def\p1{{\mathbb{P}^1_\zbb}}

\newtheorem{Theorem}{\quad Theorem}[section]

\newtheorem{Lemma}[Theorem]{\quad Lemma}

\begin{document}

\title{About Brezis-Merle Problem with Holderian condition: the case of three blow-up points.}

\author{Samy Skander Bahoura}

\address{ Departement de Mathematiques, Universite Pierre et Marie Curie, 2 place Jussieu, 75005, Paris, France.}
              
\email{samybahoura@yahoo.fr} 


\date{}

\maketitle

\begin{abstract}

We consider the following problem on open set  $ \Omega $ of $ {\mathbb R}^2 $:

\begin{displaymath} \left \{ \begin {split} 
      -\Delta u_i & = V_i e^{u_i} \,\, &\text{in} \,\, &\Omega \subset {\mathbb R}^3, \\
                  u_i  & = 0  \,\,  & \text{on} \,\, &\partial \Omega.                          
\end {split}\right.
\end{displaymath}

We assume that :

$$ \int_{\Omega} e^{u_i} dy  \leq C, $$

and,

$$ 0 \leq V_i \leq b  < + \infty $$

On the other hand, if  we assume that $ V_i  $ $ s- $holderian with $ 1/2< s \leq 1$, and,  

$$ \int_{\Omega}V_i e^{u_i} dy \leq 32\pi-\epsilon , \,\, \epsilon >0 $$
 
 then we have a compactness result, namely:

$$ \sup_{\Omega} u_i \leq c=c(b, C, A, s, \Omega). $$

where $ A $ is the holderian constant of $ V_i $.

\end{abstract}

\section{Introduction and Main Results} 

We set $ \Delta = \partial_{11} + \partial_{22} $  on open set $ \Omega $ of $ {\mathbb R}^2 $ with a smooth boundary.

\bigskip

We consider the following problem on $ \Omega \subset {\mathbb R}^2 $:

\begin{displaymath} (P) \left \{ \begin {split} 
      -\Delta u_i & = V_i e^{u_i} \,\, \\
                  u_i  & = 0  \,\,             
\end {split} 
\begin {split}
     &\text{in} \!\!&&\Omega \subset {\mathbb R}^3, \\
    & \text{in} \!\!&&\partial \Omega.               
\end {split}\right.
\end{displaymath}

We assume that,

$$ \int_{\Omega} e^{u_i} dy  \leq C, $$

and,

$$ 0 \leq V_i \leq b  < + \infty $$

The previous equation is called, the Prescribed Scalar Curvature equation, in
relation with conformal change of metrics. The function $ V_i $ is the
prescribed curvature.

\bigskip

Here, we try to find some a priori estimates for sequences of the
previous problem.

\smallskip

Equations of this type were studied by many authors, see [2, 3, 4, 6, 7, 8, 10, 11, 12, 13, 16]. We can
see in [4], different results for the solutions of those type of
equations with or without boundaries conditions and, with minimal
conditions on $ V $, for example we suppose $ V_i \geq 0 $ and  $ V_i
\in L^p(\Omega) $ or $ V_ie^{u_i} \in L^p(\Omega) $ with $ p \in [1,
+\infty] $. 

\medskip

Among other results, we  can see in [4], the following important Theorem,

\medskip

{\bf Theorem A} {\it (Brezis-Merle [4])}.{\it If $ (u_i)_i $ and $ (V_i)_i $ are two sequences of functions relatively to the previous problem $ (P) $ with, $ 0 < a \leq V_i \leq b < + \infty $, and without the boundary condition,  then, for all compact set $ K $ of $ \Omega $,

$$ \sup_K u_i \leq c = c(a, b, m, K, \Omega) \,\,\, {\rm if } \,\,\, \inf_{\Omega} u_i \geq m. $$}

A simple consequence of this theorem is that, if we assume $ u_i = 0 $ on $ \partial \Omega $ then, the sequence $ (u_i)_i $ is locally uniformly bounded. We can find in [4] an interior estimate if we assume $ a=0 $, but we need an assumption on the integral of $ e^{u_i} $. We have in [4]:

\smallskip

{\bf Theorem B} {\it (Brezis-Merle [4])}.{\it If $ (u_i)_i $ and $ (V_i)_i $ are two sequences of functions relatively to the previous problem $ (P) $ with, $ 0 \leq V_i \leq b < + \infty $, and,

$$ \int_{\Omega} e^{u_i} dy  \leq C, $$

then, for all compact set $ K $ of $ \Omega $,

$$ \sup_K u_i \leq c = c(b, C, K, \Omega). $$}

If, we assume $ V $ with more regularity, we can have another type of estimates, $ \sup + \inf $. It was proved, by Shafrir, see [13], that, if $ (u_i)_i, (V_i)_i $ are two sequences of functions solutions of the previous equation without assumption on the boundary and, $ 0 < a \leq V_i \leq b < + \infty $, then we have the following interior estimate:

$$ C\left (\dfrac{a}{b} \right ) \sup_K u_i + \inf_{\Omega} u_i \leq c=c(a, b, K, \Omega). $$

We can see in [7], an explicit value of $ C\left (\dfrac{a}{b}\right ) =\sqrt {\dfrac{a}{b}} $. In his proof, Shafrir has used the Stokes formula and an isoperimetric inequality. For Chen-Lin, they have used the blow-up analysis combined with some geometric type inequality for the integral curvature.

\bigskip

Now, if we suppose $ (V_i)_i $ uniformly Lipschitzian with $ A $ the
Lipschitz constant, then, $ C(a/b)=1 $ and $ c=c(a, b, A, K, \Omega)
$, see Br\'ezis-Li-Shafrir [3]. This result was extended for
H\"olderian sequences $ (V_i)_i $ by Chen-Lin, see  [7]. Also, we
can see in [10], an extension of the Brezis-Li-Shafrir to compact
Riemann surface without boundary. We can see in [11] explicit form,
($ 8 \pi m, m\in {\mathbb N}^* $ exactly), for the numbers in front of
the Dirac masses, when the solutions blow-up. Here, the notion of isolated blow-up point is used. Also, we can see in [16] refined estimates near the isolated blow-up points and the  bubbling behavior  of the blow-up sequences.

\bigskip

In [4], Brezis and Merle proposed the following Problem:

\smallskip

{\bf Problem} {\it (Brezis-Merle [4])}.{\it If $ (u_i)_i $ and $ (V_i)_i $ are two sequences of functions relatively to the previous problem $ (P) $ with, 

$$ 0 \leq V_i \to V \,\, {\rm in } \,\,  C^0(\Omega). $$

$$ \int_{\Omega} e^{u_i} dy  \leq C, $$

Is it possible to prove that:

$$ \sup_{\Omega} u_i \leq c = c(C, V, \Omega) \,\, ? $$}

Here, we assume more regularity on $ V_i $, we suppose that $ V_i \geq 0 $ is $ C^s $ ($ s$-holderian) $ (1/2 < s \leq 1)$ . We give the answer where $ b C < 32 \pi $.

\bigskip

In the similar way, we have in dimension $ n \geq 3 $, with different methods, some a priori estimates of the type $ \sup \times \inf $ for equation of the type:

$$ -\Delta u+ \dfrac{n-2}{4(n-1)}  R_g(x)u = V(x) u^{(n+2)/(n-2)} \,\,\, {\rm on } \,\, M. $$

where $ R_g $ is the scalar curvature of a riemannian manifold $ M $, and $ V $ is a function. The operator $  \Delta=\nabla^i(\nabla_i) $ is the Laplace-Beltrami operator on $ M $.

\bigskip

When $ V \equiv 1 $ and $ M $ compact, the previous equation is the Yamabe equation. T. Aubin and R. Scheon solved the Yamabe problem, see for example [1]. If $ V $ is not a constant function, the previous equation is called a prescribing curvature equation, we have many existence results see also [1], for a detailed summary. 

\bigskip

We can see in [3], [5], [9], some results for elliptic equations of this type,  and,  some application of the method of moving-plane  to obtain uniform estimates and estimates of type $ \sup \times \inf $. See also, [14], [15] for other $ \sup + \inf $ inequalities.

In [5], we have a classification result for singular and non-singular solution of the Yamabe equation on open set of $ {\mathbb R}^n $ and on $ {\mathbb R}^n $. The method used is of moving-plane and some other estimates. 

In [9], we have a basic description of the method of moving-plane, and, in [3], we have an application of this method, namely; inequality of type $  \sup + \inf  $ on a bounded domain of $ {\mathbb R}^2 $.

\bigskip

Returning to ourk wrok,  we give a compactness result for the Brezis-Merle Problem when the energy is less than $ 32 \pi - \epsilon $, $ \epsilon >0$. In fact, we extend the result of the author, see [2]. We argue by contradiction and we use some asymptotic estimates for the blow-up functions. Also, we use a term of the Pohozaev identity to conclude to a contradiction.


\bigskip

Our main result is:

\bigskip

{\bf Theorem}. {\it Assume that, $ V_i  $ is uniformly $ s- $holderian with $ 1/2 < s \leq 1$, and,

$$ \int_{B_1(0)} V_i e^{u_i}  dy  \leq 32\pi-\epsilon, \,\, \epsilon>0,$$

then we have:

$$ \sup_{\Omega} u_i \leq c=c(b, C, A, s, \Omega). $$

for solutions of the problem $ (P) $, here $ A $ is the holderian constant of $ V_i $.}

\section{Proof of the result:} 

\underbar {Proof of the theorem:}

\bigskip

Without loss of generality, we can assume that $ \Omega= B_1(0) $ the unit ball centered on the origin.

\bigskip

Here, $ G $ is the Green function of the Laplacian with Dirichlet condition on $ B_1(0) $. We have (in complex notation):
 
 $$ G(x, y)=\dfrac{1}{2 \pi} \log \dfrac{|1-\bar xy|}{|x-y|}, $$ 
 
we can write:

$$ u_i(x)=\int_{B_1(0)} G(x,y) V_i(y) e^{u_i(y)} dy , $$
 
We assume that we are in the case of one blow-up point. Following the notation of a previous paper, see [2], we have:

$$ \max_{\Omega} u_i=u_i(x_i) \to + \infty, $$

$$ \delta_i = d(x_i, \partial \Omega) \to 0, $$

 for $ \epsilon >0 $ small enough, and $ |y| =  \epsilon $,
 
  $$ u_i(x_i+\delta_i y) \leq C_{\epsilon}, $$
  
 and, the $ \sup +\inf $ inequality gives:
 
 $$ u_i(x_i)+4\log \delta_i \leq C. $$
 
 Also, we have the following estimates which imply the smallness for a term of the Pohozaev identity:
 
 $$  ||\nabla u_i ||_{L^q(B(x_i, \delta_i \epsilon' ))} =o(1). \,\, \forall \,\, 1 \leq q < 2. $$

 We have; because $ V_i $ is s-holderian with $ 1/2 < s \leq 1 $, the following term of the Pohozaev identity tends to $ 0 $

 $$  J_i= \int_{B(x_i, \delta_i \epsilon')} <x_1^i |\nabla (u_i-u)> (V_i-V_i(x_i)) e^{u_i} dy = o(1). $$

Now, we set:

$$ r_i = e^{-u_i(x_i)/2}, $$

we write, for $ |\theta| \leq \dfrac{\delta_i\epsilon'}{r_i}  $, $ 0 < \epsilon' <\frac{1}{4} $,

$$ u_i(x_i+r_i \theta)=\int_{\Omega} G(x_i+r_i\theta,y) V_i(y) e^{u_i(y)} dx= $$

$$ =\int_{\Omega-B(x_i, 2\delta_i\epsilon')} G(x_i,y) V_i e^{u_i(y)} dy + \int_{B(x_i, 2\delta_i\epsilon')} G(x_i+r_i\theta,y) V_ie^{u_i(y)} dy = $$ 
 
We write, $ y=x_i+r_i\tilde \theta $, with $ |\tilde \theta| \leq 2\dfrac{\delta_i}{r_i}\epsilon' $,

 $$u_i(x_i+r_i \theta) =\int_{B(0, 2\frac{\delta_i}{r_i}\epsilon')} \dfrac{1}{2 \pi} \log \dfrac{|1-(\bar x_i+ r_i \bar \theta )(x_i+r_i \tilde \theta )|}{r_i|\theta-\tilde \theta|} V_ie^{u_i(y)} r_i^2dy  + $$
 
 $$ + \int_{\Omega-B(x_i, 2\delta_i\epsilon')} G(x_i+r_i\theta,y) V_i e^{u_i(y)} dy $$

 $$ u_i(x_i)=\int_{\Omega-B(x_i, 2\delta_i\epsilon')} G(x_i,y) V_i e^{u_i(y)} dy + \int_{B(x_i, 2\delta_i\epsilon' )} G(x_i,y) V_ie^{u_i(y)} dy $$
  
 Hence,
 
 $$  u_i(x_i)=\int_{B(0, 2\frac{\delta_i}{r_i}\epsilon')} \dfrac{1}{2 \pi} \log \dfrac{|1-\bar x_i(x_i+r_i \tilde \theta )|}{r_i|\tilde \theta|} V_ie^{u_i(y)} r_i^2dy  + $$
 
 $$ + \int_{\Omega-B(x_i,2\delta_i\epsilon')} G(x_i,y) V_i e^{u_i(y)} dy  $$
   
We look to the difference,

$$ v_i(\theta)= u_i( x_i+r_i\theta)-u_i(x_i)=\int_{B(0, 2\frac{\delta_i}{r_i}\epsilon')} \dfrac{1}{2 \pi} \log \dfrac{|\tilde \theta |}{|\theta-\tilde \theta|} V_ie^{u_i(y)} r_i^2dy  +  h_1+ h_2, $$

where,

 $$ h_1(\theta) = \int_{\Omega-B(x_i, 2\delta_i\epsilon')} G(x_i+r_i\theta,y) V_i e^{u_i(y)} dy -  \int_{\Omega-B(x_i, 2\delta_i\epsilon')} G(x_i,y) V_i e^{u_i(y)} dy, $$  
 
 and,
 
 $$ h_2(\theta)= \int_{B(0, 2\delta_i\epsilon')} \dfrac{1}{2 \pi} \log \dfrac{|1-(\bar x_i+ r_i \bar \theta )y|}{|1-\bar x_i y |} V_ie^{u_i(y)} dy. $$
 
 Remark that, $ h_1 $ and $ h_2 $ are two harmonic functions, uniformly bounded.
 
 \bigskip
  
According to the maximum principle,  the harmonic function $ G(x_i+r_i\theta, .) $ on $ \Omega-B(x_i,2\delta_i\epsilon') $ take its maximum on the boundary of $ B(x_i, 2\delta_i\epsilon') $, we can compute this maximum:

$$ G(x_i+r_i\theta, y_i)=\dfrac{1}{2 \pi} \log \dfrac{|1-(\bar x_i+ r_i\bar \theta)y_i|}{|x_i+r_i\theta-y_i|} \simeq \dfrac{1}{2\pi} \log \dfrac{(|1+|x_i|)\delta_i-\delta_i(3\epsilon'+o(1))|}{\delta_i \epsilon' } \leq C_{\epsilon'} < + \infty $$

with $ y_i=x_i+ 2\delta_i \theta_i \epsilon' $,  $ |\theta_i|= 1 $, and $ |r_i\theta| \leq \delta_i\epsilon' $.

\bigskip

We can remark, for $ |\theta| \leq \dfrac{\delta_i\epsilon'}{r_i}  $, that $ v_i $ is such that:

$$ v_i = h_1+h_2+\int_{B(0, 2\frac{\delta_i}{r_i}\epsilon')} \dfrac{1}{2 \pi} \log \dfrac{|\tilde \theta |}{|\theta-\tilde \theta|} V_ie^{u_i(y)} r_i^2dy, $$

$$ v_i = h_1+h_2+\int_{B(0,2\frac{\delta_i}{r_i} \epsilon')} \dfrac{1}{2 \pi} \log \dfrac{|\tilde \theta |}{|\theta-\tilde \theta|} V_i (x_i+ r_i \tilde \theta) e^{v_i(\tilde \theta)} d\tilde \theta, $$

with $ h_1 $ and $ h_2 $, the two uniformly bounded harmonic functions.

\bigskip

{\bf Remark:} In the case of 2 or 3 blow-up points, and if we consider the half ball, we have supplemntary terms,  around the 2 other blow-up terms. Note that the Green function of the half ball is quasi-similar to the one of the unit ball and our computations are the same if we consider the half ball.

\bigskip

 We assume that, the blow-up limit is $ 0 $ and we take:
 
$$ G(x, y)=\dfrac{1}{2 \pi} \log \dfrac{|1-\bar xy|}{|x-y|}- \dfrac{1}{2 \pi} \log \dfrac{|1-xy|}{|\bar x-y|}, $$.

\underbar { \bf Asymptotic estimates and the case of one, two and three blow-up points :}

\bigskip

By the asymptotic estimates of Cheng-Lin see [8], we can see that, we have the following uniform estimates at infinity:

\begin{Lemma}

 $$ \forall \,\, \epsilon, \epsilon'>0, \,\, \exists \,\, k_{\epsilon, \epsilon'} \in {\mathbb R}_+ , \,\, i_{\epsilon, \epsilon'} \in {\mathbb N}  \,\, {\rm and} \,\,  C_{\epsilon, \epsilon'} >0 , \,\, {\rm such \, that \, for } \,\, i \geq i_{\epsilon, \epsilon' }  \,\, {\rm and} \,\, k_{\epsilon, \epsilon' } \leq  |\theta| \leq \dfrac{\delta_i\epsilon'}{r_i}  $$
 
 $$  (-4-\epsilon)\log | \theta | -C_{\epsilon, \epsilon'}  \leq v_i(\theta) \leq (-4+\epsilon)\log | \theta | + C_{\epsilon,\epsilon'}, $$
\end{Lemma} 
 
 For the proof, we consider the three following sets:

$$ A_1=\{ \tilde \theta, | \tilde \theta|\leq k_{\epsilon}\}, \,\, A_2= \{ \tilde \theta, |\theta-\tilde \theta| \leq \dfrac{|\theta|}{2}, | \tilde \theta|\geq k_{\epsilon}\}, $$

and,

$$ A_3= \{ \tilde \theta, |\theta-\tilde \theta| \geq \dfrac{|\theta|}{2}, | \tilde \theta|\geq k_{\epsilon}\}. $$

where $ k_{\epsilon} $ is such that;

$$  8\pi (1- \epsilon) \leq \int_{B(0, k_{\epsilon})}V_i (x_i+ r_i \tilde \theta) e^{v_i(\tilde \theta)} d\tilde \theta = \int_{B(x_i, k_{\epsilon} e^{-u_i(x_i)/2})} V_i e^{u_i(y)} dy \leq 8\pi (1+ \epsilon). $$
  
 In fact, if we assume that we have one blow-up point:
 
 $$  \int_{B(0,\frac{\delta_i}{2r_i})}V_i (x_i+ r_i \tilde \theta) e^{v_i(\tilde \theta)} d\tilde \theta = \int_{B(x_i,\frac{\delta_i}{2})} V_i e^{u_i(y)} dy \to 8\pi, $$
  
 To have the uniform bounds $ C_{\epsilon}>0 $, we need to bound uniformly the following quantity:

  $$ A_i= \int_{B(0,\frac{\delta_i}{2r_i})} \dfrac{1}{2 \pi} \log |\tilde \theta |V_ie^{u_i(y)} r_i^2dy = \int_{B(0,\frac{\delta_i}{2r_i})} \dfrac{1}{2 \pi} \log |\tilde \theta |V_ie^{v_i(\tilde \theta)}d\tilde \theta. $$
 
To obtain this uniform bound, we use the CC.Chen and C.S. Lin computations, see [7], to have the existence of a sequence $ l_i \to + \infty $ such that:

$$ \int_{B(0, l_i)} \dfrac{1}{2 \pi} \log |\tilde \theta |V_ie^{v_i(\tilde \theta)}d\tilde \theta \leq C, $$

and, on the other hand, the computations of YY.Li and I. Shafrir, see [11],  to have, for $ l_i \leq  |\tilde \theta| \leq \dfrac{\delta_i}{2r_i}  $ :

$$ e^{v_i(\tilde \theta)}\leq   \dfrac{C}{|\tilde \theta |^{2\beta+2} }, $$
 
for some $ 0 < \beta <1 $.

\bigskip

Finaly,

$$ A_i \leq C. $$

 \bigskip
 
Remark that, in the estimate of CC.Chen and C.S Lin, see [7], we need the assumption that $ V_i $ is $ s-$ h¬olderian with $ 0 < s \leq 1 $. 
 
 To explain more the previous lemma, we write:

$$ -2\pi v_i +2\pi h_1+2\pi h_2 = -\int_{B(0,2\frac{\delta_i}{r_i} \epsilon') \cap A_1} \dfrac{1}{2 \pi} \log \dfrac{|\tilde \theta |}{|\theta-\tilde \theta|} V_i (x_i+ r_i \tilde \theta) e^{v_i(\tilde \theta)} d\tilde \theta + , $$
 
 $$ - \int_{B(0,2\frac{\delta_i}{r_i} \epsilon') \cap A_2} \dfrac{1}{2 \pi} \log \dfrac{|\tilde \theta |}{|\theta-\tilde \theta|} V_i (x_i+ r_i \tilde \theta) e^{v_i(\tilde \theta)} d\tilde \theta + $$
 
$$ -  \int_{B(0,2\frac{\delta_i}{r_i} \epsilon') \cap A_3} \dfrac{1}{2 \pi} \log \dfrac{|\tilde \theta |}{|\theta-\tilde \theta|} V_i (x_i+ r_i \tilde \theta) e^{v_i(\tilde \theta)} d\tilde \theta = $$

$$ = -I_1- I_2- I_3. $$

For $ I_2 $, we have: $ |\theta-\tilde \theta| \leq |\tilde \theta| $, hence,

$$ -I_2 \leq 0. $$

For $ I_1 $, it is easy to see that:

$$ -I_1 \leq \log |\theta | \int_{A_1}V_i (x_i+ r_i \tilde \theta) e^{v_i(\tilde \theta)} d\tilde \theta + C, $$

with $ C $ a constant independant of $ x $ and $ i $. Here we use the estimates of Chen-Lin.

Since, $  |\theta-\tilde \theta| \leq |\tilde \theta| + |\theta | \leq |\tilde \theta|  |\theta| $ for  $ |\theta|, |\tilde \theta| \geq 1 $, we have:

$$ -I_3 \leq \log |\theta | \int_{A_3}V_i (x_i+ r_i \tilde \theta) e^{v_i(\tilde \theta)} d\tilde \theta, $$  

Thus,

$$ -2\pi v_i +2\pi h_1+2\pi h_2 \leq \log |\theta | \int_{A_1\cup A_3}V_i (x_i+ r_i \tilde \theta) e^{v_i(\tilde \theta)} d\tilde \theta + C, $$

Hence,

$$ -2\pi v_i +2\pi h_1+2\pi h_2 \leq (8\pi + \tilde \epsilon) \log |\theta | + C, $$

Thus,

$$  v_i-h_1-h_2 \geq (-4 - \epsilon) \log |\theta | - C. $$

For the rest of the proof, we use the same argument as in Cheng-Lin, see [8].

We write:

$$ v_i-h_1-h_2= \int_{B(0,2\frac{\delta_i}{r_i} \epsilon') \cap A_1} \dfrac{1}{2 \pi} \log \dfrac{|\tilde \theta |}{|\theta-\tilde \theta|} V_i (x_i+ r_i \tilde \theta) e^{v_i(\tilde \theta)} d\tilde \theta + , $$
 
 $$ + \int_{B(0,2\frac{\delta_i}{r_i} \epsilon') \cap A_2} \dfrac{1}{2 \pi} \log \dfrac{|\tilde \theta |}{|\theta-\tilde \theta|} V_i (x_i+ r_i \tilde \theta) e^{v_i(\tilde \theta)} d\tilde \theta + $$
 
$$  + \int_{B(0,2\frac{\delta_i}{r_i} \epsilon') \cap A_3} \dfrac{1}{2 \pi} \log \dfrac{|\tilde \theta |}{|\theta-\tilde \theta|} V_i (x_i+ r_i \tilde \theta) e^{v_i(\tilde \theta)} d\tilde \theta = $$

$$ = I_1 + I_2 + I_3. $$

We have:

$$ I_1 \leq - \log |\theta | \int_{A_1} \dfrac{1}{2 \pi} V_i (x_i+ r_i \tilde \theta) e^{v_i(\tilde \theta)} d\tilde \theta + C, $$

with $ C $ a constant independant of $ x $ and $ i $. Here we use the estimates of Chen-Lin.

For $ I_3 $, we have:

$$ I_3 \leq 1. $$

For $ I_2 $, we have:

$$ I_2 \leq   \dfrac{1}{2 \pi} \int_{ \{ |\tilde \theta- \theta | \leq {|\theta |}^{-\sigma} \} } \log \dfrac{1}{|\theta-\tilde \theta|} V_i (x_i+ r_i \tilde \theta) e^{v_i(\tilde \theta)} d\tilde \theta + \dfrac{\epsilon}{2} \log |\theta |, $$  

Hence,

$$  v_i-h_1-h_2 \leq (-4+ \epsilon ) \log |\theta | + \dfrac{1}{2 \pi} \int_{ \{ |\tilde \theta- \theta | \leq {|\theta |}^{-\sigma} \} } \log \dfrac{1}{|\theta-\tilde \theta|} V_i (x_i+ r_i \tilde \theta) e^{v_i(\tilde \theta)} d\tilde \theta. $$

As in [8], we can prove that, ($ h_1$ and $ h_2 $ are harmonic and satisfy the mean value theorem):

$$ v_i-h_1-h_2- \int_{ \{ |\tilde \theta- \theta | =r= {|\theta |}^{-\sigma} \} } (v_i-h_1-h_2) = \dfrac{1}{2 \pi} \int_{ B_r(x) } \log \dfrac{r}{|\theta-\tilde \theta|} V_i (x_i+ r_i \tilde \theta) e^{v_i(\tilde \theta)} d\tilde \theta $$

$$ \int_{ \{ |\tilde \theta- \theta | =r= {|\theta |}^{-\sigma} \} } (v_i-h_1-h_2) \leq (-4+ \epsilon) \log |\theta |. $$

As in the proof of the theorem 1.1 of [8], we use the Brezis-Merle estimate and the two previous estimates to prove that for $ \theta $ large enough, we have:

$$ v_i-h_1-h_2  \leq (-4+ \epsilon ) \log |\theta | + C. $$

To see this : (We write $ v_i -h_1-h_2= k+q $, with $ q $ harmonic with the same boundary  value as $ v_i $, we use Brezis-Merle estimate). Note that, $ h_1 $ and $ h_2 $ are uniformly bounded. Let $ \Omega =B_r(\theta)  $, where $ r=2|\theta |^{-\sigma} $ we have:

$$  \begin {cases}

     -\Delta k = V_i e^{k+q} & \text{in} \,\,\Omega, \\
     
 \,\, \,\,  \quad k  = 0   & \text{on} \,\, \partial \Omega.
                  
              \end {cases} $$
By the Brezis-Merle estimate:

$$  \int_{ \Omega } e^{2k} \leq C_1 |\theta |^{-2\sigma}. $$

We use the fact that $ q $ is harmonic to have:

$$ q(\theta) \leq C q(0)+ (C-1) (-\min_{ \Omega} q^-). $$

By the previous computations we have:

$$ \min_{ \Omega} q^- = \min_{\partial  \Omega} q^-= \min_{\partial  \Omega} (v_i-h_1-h_2)^- \geq (-4 - \epsilon) \log |\theta | - C, $$

and by the previous mean value estimate, we have:

$$ q(0)= \int_{ \{ |\tilde \theta- \theta | =r= {|\theta |}^{-\sigma} \} } (v_i-h_1-h_2) \leq (-4+ \epsilon) \log |\theta |. $$

Thus,

$$  q(\theta) \leq C \log |\theta |. $$

Here $ C $ is a constant independant of $ i $ and $ \sigma $.

We have by the same computations as in the proof of the theorem 1.1 of [8] to conclude that:

$$ \int_{ \{ |\tilde \theta- \theta | \leq {|\theta |}^{-\sigma} \} } e^{2 v_i} \leq |\theta |^{-2\sigma + 2 C }, $$

and by Cauchy-Schwarz inequality, we have:

 $$ \left ( \dfrac{1}{2 \pi} \int_{ \{ |\tilde \theta- \theta | \leq {|\theta |}^{-\sigma} \} } \log \dfrac{1}{|\theta-\tilde \theta|} V_i (x_i+ r_i \tilde \theta) e^{v_i(\tilde \theta)} d\tilde \theta \right )^2  \leq C, $$
 
 and that, for $ \theta $ and $ \sigma $ large enough:
  
 $$ v_i-h_1-h_2  \leq (-4+ \epsilon ) \log |\theta | + C. $$

\bigskip

Now, we extend the previous asymptotic estimates to the first derivatives:

we have, after derivation under the integral:

$$  \partial_j v_i =\partial_j h_1+\partial_j h_2+\int_{B(0, 2\frac{\delta_i}{r_i}\epsilon')} \dfrac{1}{2 \pi}  \dfrac{\theta_j-\tilde \theta_j}{|\theta-\tilde \theta|^2} V_ie^{u_i(y)} r_i^2dy, $$

In other words, we have:

$$  \partial_j v_i =\partial_j h_1+\partial_j h_2+\int_{B(0,2\frac{\delta_i}{r_i}\epsilon')} \dfrac{1}{2 \pi}  \dfrac{\theta_j-\tilde \theta_j}{|\theta-\tilde \theta|^2} V_i (x_i+ r_i \tilde \theta) e^{v_i(\tilde \theta)} d\tilde \theta, $$





We have the following lemma:

\begin{Lemma}

 $ \forall \,\, \epsilon, \epsilon'>0 \,\, \exists \,\, k_{\epsilon,\epsilon'} \in {\mathbb R}_+ , \,\, i_{\epsilon, \epsilon'} \in {\mathbb N} $,  such that, for $ i \geq i_{\epsilon, \epsilon'} $ and $ k_{\epsilon, \epsilon'} \leq  |\theta| \leq \dfrac{\delta_i\epsilon'}{r_i}  $,
 
$$  \partial_j v_i(\theta) \simeq \partial_j u_0(\theta)  \pm \dfrac{\epsilon}{|\theta|} + C\left (\dfrac{r_i}{\delta_i}\right ) ,  $$

where $ u_0 $ is the solution to:

$$ -\Delta u_0 = V(0) e^{u_0}, \,\,\, \,\, \text{in} \,\, {\mathbb R}^2. $$

\end{Lemma}

For the proof, we consider the three following sets:

$$ A_1=\{ \tilde \theta, | \tilde \theta|\leq k_{\epsilon}\}, \,\, A_2= \{ \tilde \theta, |\theta-\tilde \theta| \leq \dfrac{|\theta|}{2}, | \tilde \theta|\geq k_{\epsilon}\}, $$

and,

$$ A_3= \{ \tilde \theta, |\theta-\tilde \theta| \geq \dfrac{|\theta|}{2}, | \tilde \theta|\geq k_{\epsilon}\}. $$

where $ k_{\epsilon} $ is such that;

$$  8\pi (1- \epsilon) \leq \int_{B(0, k_{\epsilon})}V_i (x_i+ r_i \tilde \theta) e^{v_i(\tilde \theta)} d\tilde \theta = \int_{B(x_i, k_{\epsilon} e^{-u_i(x_i)/2})} V_i e^{u_i(y)} dy \leq 8\pi (1+ \epsilon). $$

{\bf Remark 1:} In the case of 2 or 3 blow-up points, and if we consider the half ball, we have supplemntary terms,  around the 2 other blow-up terms. Note that the Green function of the half ball is quasi-similar to the one of the unit ball. In the case of 3 blow-up points, we have the following supplementary term ( $ x_i $ is the principal blow-up point and $ y_i $ and $ t_i $ the 2 other blow-up points):

$$ C_1\left (\dfrac{r_i}{d(x_i, y_i)} \right ) + C_2 \left (\dfrac{r_i}{d(x_i,t_i)} \right ). $$

\bigskip

 We assume that, the blow-up limit is $ 0 $ and we take:
 
$$ G(x, y)=\dfrac{1}{2 \pi} \log \dfrac{|1-\bar xy|}{|x-y|}- \dfrac{1}{2 \pi} \log \dfrac{|1-xy|}{|\bar x-y|}, $$.









In the previous computations, we have considered the unit ball, but by a conformal transformation , we can have the same estimates on the half ball, with a coefficient of the conformal transformation. We can assume the estimates on the half ball.

\bigskip

Now, we consider the following term of the Pohozaev identity

$$  J_i= \int_{B(x_i, \delta_i \epsilon')} <x_1^i |\nabla (u_i-u)> (V_i-V_i(x_i)) e^{u_i} dy, $$

We want to show that this term tends to $ 0 $ as  $ i $ tends to infinity. We can reduce the problem, after integration by parts, to the following integral:

$$  J'_i= \delta_i \int_{B(x_i, \delta_i \epsilon')}  \partial_1 u_iV_i e^{u_i} dy = \delta_i \int_{B(x_i, \delta_i \epsilon')}  \partial_1 u_i(-\Delta u_i) $$

But, if we take $ y=x_i+r_i \theta $, with, $ |\theta| \leq \dfrac{\delta_i\epsilon'}{r_i} $,  we have:

$$ J'_i= \dfrac{\delta_i}{r_i} \int_{B(0, \frac{\delta_i}{r_i} \epsilon')}  \partial_1 v_i(-\Delta v_i) = $$

$$ = \dfrac{\delta_i}{r_i} \int_{\partial B(0, \frac{\delta_i}{r_i} \epsilon')} \left ( \frac{(\partial_1 v_i)^2}{2} \nu_1-\frac{(\partial_2 v_i)^2}{2} \nu_1+ (\partial_1 v_i)(\partial_2 v_i) \nu_2 \right ) d\sigma_i, $$

Thus, if we use the uniform asymptotic estimates, we can see that, we reduce the computation to the Pohozaev identity for the limit blow-up function (which equal to $ 0 $),  plus terms in $ {\epsilon}{|\theta|} $ and $|\theta| $. First, we tend $ i $ to infinity, after $ \epsilon $ to  $ 0$ and finaly , we tend $ \epsilon' $ to $ 0 $ .

With this method we can have a compactness result for 3 blow-ups points. First, we can see the case of 3 exteriors blow-up points, then by the previous formulation we have a compactness result,  it is the case for one of the following cases ( if we set $ \delta_i, \delta_i'  $ and $ \delta_i'' $ for the radii of each exterior blow-up) :

$$  \dfrac{d(x_i,y_i)}{\delta_i}\to + \infty \,\,\, {\rm and } \,\,\,\dfrac{d(x_i,t_i)}{\delta_i}\to + \infty, $$

or,

$$  \dfrac{d(y_i,x_i)}{\delta_i'}\to + \infty \,\,\, {\rm and } \,\,\,\dfrac{d(y_i,t_i)}{\delta_i'}\to + \infty, $$

or,

$$  \dfrac{d(t_i,x_i)}{\delta_i''}\to + \infty \,\,\, {\rm and } \,\,\,\dfrac{d(t_i,y_i)}{\delta_i''}\to + \infty, $$

or, 

the case when the distance to two exterior blow-up points is of order the radii. In this last case, we divide the region in 3 parts and use the Pohozaev identity directly. In fact, we are reduced to the case of two blow-up points.

In fact, in the case of 3 exterior blow-up points. By the previous formulation around each exterior blow-up point we look to the one of the 3 first cases.  For example, assume the first case. Then we work around the first blow-up. In fact we have, for 3 blow-up points :

\begin{Lemma} $ \forall \,\, \epsilon>0, \,\epsilon'>0 \,\, \exists \,\, k_{\epsilon,\epsilon'} \in {\mathbb R}_+ , \,\, i_{\epsilon, \epsilon'} \in {\mathbb N} $ and $ C_{\epsilon, \epsilon'} >0 $, such that, for $ i \geq i_{\epsilon,\epsilon'} $ and $ k_{\epsilon, \epsilon'} \leq  |\theta| \leq \dfrac{\delta_i\epsilon'}{r_i}  $,
 
 $$  (-4-\epsilon)\log | \theta | -C_{\epsilon, \epsilon'}  \leq v_i(\theta) \leq (-4+\epsilon)\log | \theta | + C_{\epsilon, \epsilon'}, $$
 
 and,
 
$$ \partial_j v_i \simeq \partial_j u_0(\theta) \pm \dfrac{\epsilon}{|\theta|} + C\left (\dfrac{r_i}{\delta_i} \right )^2|\theta| + m\times \left (\dfrac{r_i}{\delta_i} \right ) + C_1\left (\dfrac{r_i}{d(x_i, y_i)} \right ) + C_2 \left (\dfrac{r_i}{d(x_i,t_i)} \right )$$

\end{Lemma}

\underbar {Proof of the compactness :}
\bigskip

By using the lemma we can see that we have the compactness, because:

\bigskip

(to understand this, it is sufficient to do the computations for the half ball directly by using the Green function of the half ball directly).

 We have after using the  previous term of the Pohozaev identity:
 
 $$ o(1)=J'_i= m'+ C_1o(1)+C_2o(1), $$
 
 $$ 0=  \lim_{\epsilon' } \lim_{\epsilon} \lim_i J_i'= m', $$

 which contradict the fact that $ m'>0 $.

\bigskip

\underbar {Proof of the second estimate of the lemma:}

\bigskip

 For the proof, we consider the three following sets:

$$ A_1=\{ \tilde \theta, | \tilde \theta|\leq k_{\epsilon}\}, \,\, A_2= \{ \tilde \theta, |\theta-\tilde \theta| \leq \dfrac{|\theta|}{2}, | \tilde \theta|\geq k_{\epsilon}\}, $$

and,

$$ A_3= \{ \tilde \theta, |\theta-\tilde \theta| \geq \dfrac{|\theta|}{2}, | \tilde \theta|\geq k_{\epsilon}\}. $$

where $ k_{\epsilon} $ (large enough), is such that;

$$  8\pi (1- \epsilon) \leq \int_{B(0, k_{\epsilon})}V_i (x_i+ r_i \tilde \theta) e^{v_i(\tilde \theta)} d\tilde \theta = \int_{B(x_i, k_{\epsilon} e^{-u_i(x_i)/2})} V_i e^{u_i(y)} dy \leq 8\pi (1+ \epsilon). $$

We write:

 $$  \partial_j v_i -\partial_j h_1-\partial_j h_2=\int_{B(0,2\frac{\delta_i}{r_i}\epsilon')} \dfrac{1}{2 \pi}  \dfrac{\theta_j-\tilde \theta_j}{|\theta-\tilde \theta|^2} V_i (x_i+ r_i \tilde \theta) e^{v_i(\tilde \theta)} d\tilde \theta, $$

 $$  \partial_j v_i -\partial_j h_1-\partial_j h_2=\int_{A_1} \dfrac{1}{2 \pi} \dfrac{\theta_j-\tilde \theta_j}{|\theta-\tilde \theta|^2}V_i (x_i+ r_i \tilde \theta) e^{v_i(\tilde \theta)} d\tilde \theta + \int_{A_2 \cup A_3} \dfrac{1}{2 \pi} \dfrac{\theta_j-\tilde \theta_j}{|\theta-\tilde \theta|^2} V_i (x_i+ r_i \tilde \theta) e^{v_i(\tilde \theta)} d\tilde \theta $$

Using the estimates of $ v_i $, we obtain:

$$  \partial_j v_i -\partial_j h_1-\partial_j h_2=\dfrac{o(1)}{|\theta |}+ \int_{A_1}  \dfrac{1}{2 \pi} \dfrac{\theta_j-\tilde \theta_j}{|\theta-\tilde \theta|^2}V_0 e^{u_0(\tilde \theta)} d\tilde \theta + \int_{A_2 \cup A_3} \dfrac{1}{2 \pi} \dfrac{\theta_j-\tilde \theta_j}{|\theta-\tilde \theta|^2} V_i (x_i+ r_i \tilde \theta) e^{v_i(\tilde \theta)} d\tilde \theta $$

Thus,

$$  \partial_j v_i -\partial_j h_1-\partial_j h_2= \partial_j u_0+ \dfrac{o(1)}{|\theta |}+ \int_{A_2 \cup A_3} \dfrac{1}{2 \pi} \dfrac{\theta_j-\tilde \theta_j}{|\theta-\tilde \theta|^2} V_i (x_i+ r_i \tilde \theta) e^{v_i(\tilde \theta)} d\tilde \theta $$

Finaly,

$$  \partial_j v_i -\partial_j h_1-\partial_j h_2- \partial_j u_0 = \dfrac{o(1)}{|\theta |} + \int_{A_2 \cup A_3} \dfrac{1}{2 \pi} \dfrac{\theta_j-\tilde \theta_j}{|\theta-\tilde \theta|^2} V_i (x_i+ r_i \tilde \theta) e^{v_i(\tilde \theta)} d\tilde \theta. $$

For $ A_2 $ and $ A_3 $, we have:

$$  | \int_{A_3} \dfrac{1}{2 \pi} \dfrac{\theta_j-\tilde \theta_j}{|\theta-\tilde \theta|^2} V_i (x_i+ r_i \tilde \theta) e^{v_i(\tilde \theta)} d\tilde \theta | \leq \dfrac{1}{2 \pi} \int_{ \{ | \tilde \theta|\geq k_{\epsilon}\} } \dfrac{1}{|\theta |} V_i (x_i+ r_i \tilde \theta) e^{v_i(\tilde \theta)} d\tilde \theta  \leq  \dfrac{\epsilon}{|\theta |}, $$

because,

$$ \int_{ \{ | \tilde \theta|\geq k_{\epsilon}\} }V_i (x_i+ r_i \tilde \theta) e^{v_i(\tilde \theta)} d\tilde \theta  \to 0, $$

for $ k_{\epsilon} $ large enough.

For $  \theta \in A_2 $, $ |\tilde \theta| \geq \dfrac{|\theta|}{2} $ and,   $ |\tilde \theta| \geq \dfrac{|\theta- \tilde \theta|}{2} $, and we use the estimate of $ v_i $ to have:

$$  | \int_{A_2} \dfrac{1}{2 \pi} \dfrac{\theta_j-\tilde \theta_j}{|\theta-\tilde \theta|^2} V_i (x_i+ r_i \tilde \theta) e^{v_i(\tilde \theta)} d\tilde \theta |   \leq  \dfrac{C|\theta |}{|\theta |^{4 - \epsilon_0}} \leq  \dfrac{\epsilon}{|\theta |}, $$
 
for $ \theta $ large enough and $\epsilon_0 $ small enough.

\bigskip

Finaly, we have:

 $$  |\partial_j v_i -\partial_j h_1-\partial_j h_2-\partial_j u_0| \leq \dfrac{\epsilon}{|\theta |}, $$

for $ \theta $ large enough.

\bigskip

Now, it is easy to see from the definition of $ h_1 $ and $ h_2 $ that:

$$| \partial_j h_1-\partial_j h_2 -m \dfrac{r_i}{\delta_i} |\leq C_1\left (\dfrac{r_i}{d(x_i, y_i)} \right ) + C_2 \left (\dfrac{r_i}{d(x_i,t_i)} \right )$$

Thus,

$$ |\partial_j v_i -\partial_j u_0-m\dfrac{r_i}{\delta_i} | \leq C_1\left (\dfrac{r_i}{d(x_i, y_i)} \right ) + C_2 \left (\dfrac{r_i}{d(x_i,t_i)} \right )$$

for $ \theta $ large enough.

\end{document}